\input amstex

\documentstyle{amsppt}

\def\Pic{{\text{\rm Pic}}}
\def\dim{{\text{\rm dim}}}

\def\deg{{\text{\rm deg}}}
\def\PX{{\text {\rm Pic}}^0(X)}
\def\DI{{\text {\rm dim}}   }
\def\DX{{\text {\rm dim}}(X)}
\def\DY{{\text {\rm dim}}(Y)}

\def\ALB{{\text{\rm Alb}}(X)}
\def\alb{\text{\rm alb}}
\def\Alb{\text{\rm Alb}}
\def\lra{\longrightarrow}
\def\ot{{\otimes}}
\def\OO{{\Cal{O}}}
\def\PP{{\Cal{P}}}
\def\PPP{{\Bbb{P}}}
\def\ox{{\omega _X}}
\def\FF{{\Cal{F}}}

\def\QQ{{\Bbb{Q}}}
\def\PXT{{\text {\rm Pic}}^{\tau }(X)}
\def\PXA{{\text {\rm Pic}}^{\tau }(A)}
\def\CC{{\Bbb{C}}}
\def\ZZ{{\Bbb{Z}}}

\def\XX{{ {\Cal     X}}}

\topmatter

\title {Pluricanonical maps of
varieties of maximal Albanese dimension}\endtitle
\rightheadtext{Pluricanonical maps}
\author Jungkai A. Chen, Christopher D. Hacon   \endauthor

\address Jungkai Alfred Chen, Department of Mathematics, National Chung
Cheng University, Ming 
Hsiung, Chia Yi, 621, Taiwan \endaddress

\email  jkchen\@math.ccu.edu.tw  \endemail

\address Christopher Derek Hacon,  Department of Mathematics, University 
of Utah, 155 South 1400 East JWB 233, Salt
Lake City, UT 84112, USA 
\endaddress

\email   chhacon\@math.utah.edu   \endemail

\footnote""{The first author was partially supported by
National Science Council, Taiwan (NSC-89-2115-M-194-016)}

\abstract Let $X$ be a smooth complex projective
algebraic variety of maximal
Albanese dimension. We give a characterization of $\kappa (X)$ in terms
of the set $V^0(X,\omega_{X} )$ 
$:=\{ P\in {\text{\rm Pic}}^0(X)|h^0(X,  \omega_X \otimes P)  \ne 0\}$.
An immediate consequence of this is that the Kodaira dimension $\kappa (X)$
is invariant under smooth deformations.
We then study the pluricanonical maps $\varphi _m:X \dashrightarrow 
\Bbb{P} (H^0(X,mK_X))$. 
We prove that if $X$ is of general type, 
$\varphi _m$ is generically finite for $m\geq 5$ 
and birational for $m\geq 5 \text{\rm dim} (X) +1$. 
More generally, we show that for $m\geq 6$ 
the image of $\varphi _m$ is of dimension equal to $\kappa (X)$ 
and for $m\geq 6\kappa (X)+2$, 
$\varphi _m$ is the stable canonical map.

\endabstract

\date April 27, 2000 \enddate

\endtopmatter

\document

\heading Introduction \endheading
Generic vanishing theorems have recently proven to be a very effective
tool in the study of the geometry of irregular varieties.
In this paper we will show how these techniques can be used to answer a
series of natural questions about varieties of maximal Albanese dimension.

In what follows, $X$ will denote a smooth complex projective algebraic variety,
$\Alb (X)$ denotes the Albanese variety of $X$ and
$\alb _X
:X\lra \Alb (X)$ the corresponding Albanese map. We will assume that
$X$ is of maximal Albanese dimension and hence $\alb _X$ is a generically
finite morphism.
In particular $q(X)\geq \dim (X)$.
Koll\'ar made the following

\proclaim{Conjecture \cite {Ko3, (17.9.3)}}
If $X$ is of general type and maximal Albanese dimension.
Then $\chi (X, \ox)>0$.
\endproclaim
As an immediate consequence of the generic vanishing theorems of Green
and Lazarsfeld (cf. Theorem
1.3.1), one sees that this is equivalent to $h^0(X,
\omega _X \ot P)>0$
for all $P\in \Pic ^0(X)$.
In \cite {EL1, Example 1.13}, Ein and Lazarsfeld produce a counterexample
to this conjecture.
However, we will show that a similar result does hold.
\proclaim{Theorem 1}
The translates through the origin of the irreducible components of
$V^0(X,\omega _X):=\{ P\in \Pic ^0(X)|h^0(X,  \ox \ot P)  \ne 0\}$ 
generate a subvariety of $\Pic ^0(X)$ of dimension
$\kappa (X)-\dim (X)+q(X)$. In particular, if $X$ is of general type,
$V^0(X,\omega _X)$ generates $\Pic ^0(X)$.
\endproclaim
This is a generalization of a result of Ein and Lazarsfeld
\cite{EL2} (see Theorem 1.3.2).
An immediate consequence is the following
\proclaim{Corollary 2} Let $X$ be of maximal Albanese dimension.
Then the Kodaira dimension of $X$ is invariant under smooth deformations.
\endproclaim

We next turn our attention to the study of the linear series
$|mK_X|$ and the corresponding rational maps $\varphi _m:=\varphi
_{|mK_X|}$. In \cite{Ko3}, Koll\'{a}r shows
that if $X$  is  smooth proper variety
of general type with generically large algebraic fundamental group, then
the sections of $|mK_X|$
define a birational map for $m\geq 10 ^{\DX }$. Following \cite{AS}
it is reasonable to expect a quadratic bound in $\DX$.
\proclaim {Theorem 3}If $X$ is of maximal Albanese and Kodaira dimensions.
Then
$\varphi _5$ is generically finite, and $\varphi _{5\dim (X)+1}$ is
birational.
\endproclaim
Our proof is based on the techniques of Koll\'ar \cite{Ko3}.
However, technical difficulties imply that when $X$ is not of
general type, the bounds are somewhat weaker.
\proclaim{Theorem 4}
If $X$ is of maximal Albanese dimension. Then $\varphi _6(X)$ has
dimension equal to $\kappa (X)$ and $\varphi _{6 \kappa(X) +2}$ is the
stable canonical map.
\endproclaim
We do not expect these bounds to be optimal. It would be interesting
to produce examples of varieties of general type and maximal Albanese
dimension for which $\varphi _m$ is not birational for $m\geq 3$.

\bigskip
{\bf Acknowledgment.} We would like to thank 
P. Belkale, H. Clemens, Y. Kawamata,  
J. Koll\'ar, R. Lazarsfeld and I.-H. Tsai for valuable conversations.

\heading Conventions and Notations   \endheading

(0.1) Throughout this paper, we work
over the field of complex numbers $\CC$.

(0.2) For $D_1,D_2$ $\QQ$-divisors on a variety $X$, we
write $D_1 \prec D_2$ if $D_2-D_1$ is effective,
and $D_1 \equiv D_2$ if $D_1$ and $D_2$ are numerically equivalent.

(0.3) $|D|$ will denote the linear series associated to the divisor
$D$, and $Bs|D|$ denotes the base locus of $|D|$.

(0.4) For a real number $a$, let $\lfloor a \rfloor$ be the
largest integer $\le a$
and $\lceil a \rceil$ be the smallest integer $\ge a$.
For a  $\QQ$-divisor $D=\sum a_i D_i$, let
$\lfloor D \rfloor =\sum \lfloor a_i \rfloor D_i$ and
$\lceil D \rceil =\sum \lceil a_i \rceil D_i$.

(0.5) Let $\FF$ be a coherent sheaf on $X$,
then $h^i(X,\FF)$ denotes the complex dimension of $H^i(X,\FF)$.
In particular,
the plurigenera $h^0(X,\ox^{ \otimes m})$ are denoted  by $P_m(X)$
and the irregularity $h^0(X,\Omega_X^1)$ is denoted by $q(X)$.

(0.6) We will denote by $\ALB$ the albanese variety of $X$,
by $\alb _X : X\lra \ALB $ the albanese morphism.
As usual $\PX $ is the abelian variety dual to $\ALB $
parametrazing all topologically trivial line bundles on $X$.
$\PXT$ will denote the set of torsion elements in $\PX$.

(0.7) Let $f: X \to Y $ be a generically finite morphism. $R_f$
denotes the ramification divisor.

(0.8) Let $f: X \to Y $ be a morphism. An irreducible divisor $E$
on $X$ is said to be $f$-exceptional if $\DI {f(E)} < \DI {f(X)}
-1$. And $E$ is $f$-vertical if $\DI {f(E)} \le \DI {f(X)} -1$.
An effective divisor is $f$-exceptional (resp. $f$-vertical) if every
component is.

(0.9) We will say that a $\QQ$-divisor $\Delta$ is klt (Kawamata log terminal)
if $\Delta $ has normal crossings support and $\lfloor \Delta \rfloor =0$.
We refer to \cite {Ko3, 10.1.5} for the general
definition of klt divisor. 

\head 1. Preliminaries \endhead
\subhead 1.1 Higher Direct Images of Dualizing Sheaves \endsubhead

For the convenience of the reader, we recall various results
which will be frequently used in what follows.

\proclaim {Theorem 1.1.1 }
Let $X$, $Y$ be projective varieties, $X$ smooth and $f:X\lra Y$ a
surjective map.
Let $M$ be a line bundle on $X$
such that $M\equiv f^* L +\Delta$, where $L$
is a $\QQ$-divisor on $Y$
and $(X,\Delta )$ is klt. Then

(a) $R^jf_*(\omega _X\ot M) $ is torsion free for $j\geq 0$.

(b) if $L$ is nef and big, then $H^i(Y,R^jf_*(\omega _X\ot M))
=0$ for $i>0$, $j\geq 0$.

(c) $R^if_* \omega _X$ is zero for $i>\DX -\DY$.

(d) $R^{\cdot }f_*\omega _X\cong \sum R^if_*\omega _X [-i]$ and
$h^p(X,\omega _X)=\sum h^i (Y, R^{p-i}f_* \omega _X )$.

(e) If $Y$ is birational to an abelian variety and $L$ is nef and big,
it follows that $H^0(Y,R^jf_*(\omega
_X\ot M)\ot P)=0$ for all $P\in \Pic^0(Y)$
if and only if $R^jf_*(\omega _X\ot M      )=0$.
\endproclaim
\demo{Proof} See \cite{Ko1}, \cite{Ko2}, \cite{Ko3}, \cite {Ka2},
\cite {EV1}.
\hfill $\square$ \enddemo
\proclaim {Theorem 1.1.2 }
Let $f:X\lra Y$ be a surjective morphism from a smooth projective variety
to a projective
variety. Let $Q\in \PXT$. Then
$$R^{\cdot }f_*(\omega _X\ot Q)\cong \sum R^if_*(\omega _X\ot Q) [-i].$$
In particular $h^p(X,\omega _X\ot Q)=
\sum h^i (Y, R^{p-i}f_* (\omega _X \ot Q))$.
\endproclaim
\demo{Proof} \cite{Ko2, Corollary 3.3}.
\hfill $\square$ \enddemo

\bigskip
\subhead 1.2 The Iitaka Fibration \endsubhead

Let $X$ be a smooth complex projective variety of maximal Albanese
dimension.
Let $a:X\lra A$ be a generically finite map such that the image of $X$
generates the abelian variety $A$.
A nonsingular representative of the Iitaka fibration of $X$ is a morphism
of smooth complex projective varieties $f':X'\lra Y$ such that
$\dim (Y)=\kappa (X)$ and $\kappa ({X'}_y)=0$, where ${X'}_y$ is a generic
geometric fiber of $f'$. Since our questions will be birational in nature,
we may always assume that $X=X'$ and $f=f'$.
Since $\kappa ({X}_y)=0$, it follows by
\cite {Ka1} that the image of the fiber ${X_y}$ is the translate of an
abelian
subvariety $K_y$ of $A$, and $a|_{X_y}:X_y\lra K_y$ is birationally
equivalent to an \'etale map. Since $A$ can contain at most countably
many abelian subvarieties, we may assume that the $K_y$
are all translates of a fixed abelian subvariety $K$ of $A$.
We will denote by $\tilde {K}$ the abelian variety birational to $X_y$.
Let $p:A\lra S:=A/K$.
Let $Z$ (resp. $W$) denote the image of $X$ in $A$ (resp. $S$).
By construction, the general fiber of $f:X\lra Y$ maps to a closed point
in $S$, therefore replacing $X$ by an appropriate birational
model, there exists a morphism $q:Y\lra S$ such that $q\circ f = p\circ
a$.
We may assume that the above maps fit in the following commutative diagram

$$
\CD
X &@>{a}>>&Z&@>{\subset}>>&A \\
@V{f}VV  & & & & & & @VV{p}V \\
Y &@>{q}>>&W&@>{\subset}>>& S.\\
\endCD
$$

\bigskip
\subhead  1.3 Cohomological Support Loci \endsubhead

Let $a:X\to A$ be a morphism from a smooth projective variety $X$ to
an abelian variety $A$.
If $\FF$ is a coherent sheaf on $X$,
then one can define the {\it cohomological support loci} by

$$V^i(X,T,\FF):=\{ P \in T \subset  \Pic ^0(A) | h^i(X, \FF \otimes a^*P)
\ne 0
\}.$$
In particular,
if $a={ {\alb }}_X :X \to \Alb (X)$,
then we simply write
$$V^i(X,\FF):=\{ P \in \Pic ^0(X) | h^i(X, \FF \otimes P) \ne 0 \}.$$

We say that $X$ has maximal Albanese dimension if $\DI {(\alb _X(X))}=\DI
(X)$.
The geometry of the loci $V^i(X,\ox )$ defined above
is governed by the following:

\proclaim{Theorem 1.3.1 (Generic vanishing)}

(a) Any irreducible component of $V^i(X,\ox )$ is a translate of
a sub-torus of $\Pic ^0(X)$
and is of codimension at least $i-(\DI (X) -\DI (\alb _X(X)))$.

(b) Let $P$ be a general point of an irreducible component
$T$ of $V^i(X,\ox)$.
Suppose that $v
\in H^1(X, \OO
_X)\cong T_{P }\Pic ^0(X)$ is not tangent to $T$. Then the sequence
$$H^{i-1}(X,\ox \otimes P) \overset {\cup v} \to {\longrightarrow}
H^{i}(X,\ox \otimes P)  \overset {\cup v} \to {\longrightarrow }
H^{i+1}(X, \ox \otimes P) $$
is exact.
If $v$ is tangent to $T$, then the maps in the above sequence vanish.

(c) If $X$ is a variety of maximal Albanese dimension, then
$$\Pic ^0(X)\supset V^0(X,\ox )\supset V^1(X,\ox)
\supset \ ...\ \supset V^n(X,\ox )=\{ \OO _X\}.$$

(d) Every irreducible component of $V^i(X,\ox)$ is
a translate of a sub-torus of $\PX$ by a torsion point.
\endproclaim
\demo {Proof} See \cite{GL1},\cite{GL2},\cite{EL1} and \cite{S}.
\hfill $\square$ \enddemo
In \cite{EL1}, Ein and Lazarsfeld illustrate various examples
in which the geometry of $X$ can
be recovered from information on the loci $V^i(X, \ox)$.
In particular,
they prove

\proclaim{Theorem 1.3.2 \cite{EL2}}
If $X$ is a variety with
maximal Albanese dimension and $\DI V^0(X,\ox)=0$. Then $X$ is
birational to an abelian variety.
\endproclaim
\proclaim{Proposition 1.3.3 \cite {EL2}} Let $a:X\lra A$ be a generically
finite map from a smooth projective variety to an abelian variety.
Let $P$ be any isolated point of $V^0(X, \Pic ^0(A),\omega _X)$. Then
$a^* P=\OO _X$.
\endproclaim

We will also need the following
\proclaim{Lemma 1.3.4 \cite{CH, Lemma 2.1}}
Let $X$ be a variety of maximal Albanese dimension. 
Fix $Q\in \PXT$. Then $h^0(X, \ox^2 \ot Q \ot f^*P)$ is constant for all
torsion $P\in \Pic ^0(Y)$.
\endproclaim

\heading 2. Kodaira dimension of
varieties of maximal Albanese dimension  \endheading
Throughout this paper, we will assume that $X$ is of maximal Albanese
dimension and hence of positive Kodaira dimension $\kappa (X)\geq 0$.
We will need the following immediate consequence of \cite {M, Theorem 2.2}

\proclaim{Lemma 2.1}
If $\FF $ is a coherent sheaf on an abelian variety $A$
such that $h^i(A,\FF \ot P)$ $=0$ for all $i
\geq 0$ and all $P\in \Pic ^0(A)$. Then $\FF =0$.
\endproclaim
We will frequently refer to the notation and results of \S 1.
Define
$$G:=ker \left( \Pic ^0(A)\lra \Pic ^0(K)\lra
\Pic ^0(\tilde {K})\right).$$
Let $\overline {G}:=G/\Pic ^0(S)$. By dimension considerations $\overline
{G}$
is a finite group, hence $G$ consists of finitely many translates of
$\Pic ^0(S)$.
Let $Q_1,...,Q_r\in G\subset \Pic ^0(A)$
be a set of torsion line bundles representing lifts of the elements of
$\overline{G}$.

We will frequently identify $\Pic^0(S)$ with its image in $\Pic^0(A)$.

\proclaim {Lemma 2.2}

(a) $V^0(X,\Pic ^0(A),\omega _X )\subset G.$

(b) For every $Q_i$
the loci $V^0(X,Q_i+\Pic ^0(S),\omega _X )$ are nonempty.

(c) If $\OO_X \not \in a^* (Q_i+\Pic ^0(S))$, then
there exists a positive dimensional component of
$V^0(X,Q_i +\Pic ^0(S),\omega _X )$.
\endproclaim
\demo {Proof}
If $H^0(X, \omega _X \ot a^*P) \ne 0$ for some $P\in \Pic ^0(A)$,
then for general $y\in Y$,  $h^0(X_y, \omega _{X_y} \ot a^*P)=
h^0(\tilde {K}, \OO _{\tilde {K}} \ot a^*P)  \ne 0 $. 
This is possible only if $P$ is in
the kernel
of $\Pic ^0(A)\lra \Pic ^0(\tilde {K}) = \Pic^0(X_y)$.

To prove (b), consider  $\pi :X\lra W \subset S$.
Assume that $V^0(X,Q+\Pic ^0(S),\omega _X )$ is empty, then
by Theorem 1.3.1.c,
$H^i(X,\omega _X $ $\ot a ^*Q\ot  \pi ^*P)=0$ for all $i\geq 0$ and all
$P\in \Pic ^0(S)$.
By Theorem 1.1.2, $H^k(W,R^j \pi _* (\omega _X \ot $
$a ^*Q)\ot P)=0$ for all $k$, and all
$P\in \Pic ^0(S)$.
Therefore, by Lemma 2.1,
$R^j \pi _* (\omega _X \ot $ $a ^*Q)=0$ for all $j$.
In particular, for a general fiber
$X_w$ of $\pi: X\lra W$,
$$h^0(X_w,\omega _{X_w} \ot a ^*Q)=0.$$
Notice that $X_w$ is a finite union of general fibers
of $X \to Y$ hence a finite union of varieties birational to $\tilde{K}$,
therefore for $Q\in G$
$$h^0(X_w,\omega _{X_w} \ot a ^*Q)=h^0(X_w,\omega _{X_w} ) > 0,$$
which is a contradiction.
We may therefore assume that
$V^0(X,Q+\Pic ^0(S),\omega _X )$ is nonempty.

The assertion (c) now follows
since by Proposition 1.3.3,
any isolated point $Q$ of $V^0(X,\Pic ^0(A),\omega _X)$ must be such that
$a^*Q=\OO _X$.
\hfill $\square$ \enddemo
\proclaim {Theorem 2.3}
Let $a:X\lra A$ be a generically finite
morphism from a smooth complex
projective variety to an abelian variety.
Then, the translates through the origin of the components of
$V^0(X,\Pic ^0(A),\omega _X )$ generate
$\Pic ^0(S)$.
\endproclaim

\demo {Proof}
Assume that translates through the origin of the components of
$V^0:=V^0(X,\Pic ^0(A),\omega _X)$ do not generate
$\Pic ^0(S)$.
Then there exists
an abelian proper subvariety $T\subset \Pic ^0(S)$
and a finite subgroup $G'$
such that $V^0 \subset T+G'$. Consider the induced morphism
$$\pi :X\lra T^*=:C,$$
which factors through $X@>{a}>> A@>{p}>> S@>>>C$.

Let $Q\in \PXA$ be any torsion element not contained in
$T+G'$. Then, $H^0(X,\omega _X $ $\ot a ^*Q\ot  \pi ^*P)=0$ 
for all $i\geq 0$ and all
$P\in \Pic ^0(C)$.
By Theorem 1.3.1.c,
$H^i(X,\ox \ot $ $ Q\ot  \pi ^*P)=0$
for all $0\leq i\leq n$ and all $P\in \Pic ^0(C)$.

Let $V\subset C$ be the image of $X$, and $\pi:X\lra V$ be the induced map.
By Theorem 1.1.2 and the projection formula, it follows that
$$H^j(V,R^i \pi _* (\omega _X \ot Q)\ot P)=0$$ for all $i,j$ and
$P\in \Pic ^0(C)$.
By Lemma 2.1, $R^i \pi _* (\omega _X \ot Q)=0$ for all $i$.

If $\pi:X\lra V$ is generically finite, then $\pi _* (\omega _X \ot Q)$
is clearly non-zero. We may therefore assume
that $\pi$ has positive dimensional
fibers.

Let $p$ be a point in $V$ and $X_p=\pi ^{-1}(p)$
be the inverse image of $p$.
For general $p\in V$ (not depending on $Q$), we have that
$$H^i(X_p,\omega _{X_p}\ot Q)=H^i(X_p,(\omega _{X}\ot Q)\ot \OO _{X_p})=0.$$

Let $B$ be the connected component through the origin of the
kernel of $A\lra C$. The image of $X_p$ in $A$ is contained
in a translate of $B$
which we denote by $B_p$.
Since $h^0(X_p,\omega _{X_p}\ot Q)=0$ for all but finitely many $Q\in
\Pic ^0(B)$, it follows by Proposition 1.3.2 that $\kappa (X_p)=0$ and
$X_p\lra B_p$
is birationally an \'etale map of abelian varieties.

On the other hand, by the weak addition formula
$$\kappa (X_p)+\dim (V)\geq \kappa (X),$$
so $\kappa(X_p) >0$ since $\kappa(X)= \DI (S) > \DI (V)$.
This contradicts the assumption that $T$ is a proper subvariety
of $\Pic ^0(S)$.
Therefore, $T=\Pic ^0(S)$.
\hfill $\square$ \enddemo

\proclaim{Corollary 2.4} In the above hypothesis,
the dimension of the subgroup of $\Pic ^0(X)$
generated by the translates through the origin of the components
of $V^0(X,\Pic ^0(A), $ $\omega _X)$ is equal to $\kappa (X)-\dim (X)+q(X)$.
\endproclaim
\noindent {\bf Example \cite {Ko3, 17.9.5}.}
Let $p:C\lra E$ be a degree 2 map from a genus 2 curve
$C$ to an elliptic curve $E$. We may assume that $p_*\OO _C =\OO_E
\oplus L^{-1}$, for an appropriate $L\in \Pic (E)$ such that $L^{\ot 2}
=\OO _E(B)$ where $B$ is the branch locus.
Let $\tilde {F}\lra F$ be a degree 2 \'etale map of elliptic curves
such that $e_*\OO _{\tilde {F}}=\OO _F\oplus P$ with $P^{\ot 2}\cong \OO
_F$.
Define
$$X:= \tilde {F}\times C/<(i_{\tilde {F}}\times i_C)>.$$
Here $i_{\tilde {F}}$ and $i_C$ denote the involutions
on ${\tilde {F}}$ and $C$ respectively. We have that for $a:X\lra F\times
E$,
$$a_*(\omega _X)\cong (\OO _F \otimes \OO _E)\oplus (P\otimes L).$$
(Pull backs have been omitted.)
It follows that Iitaka fibration has image $E$ and 
$$V^0(X,\Pic ^0(F\times E),\omega _X )=\{ \OO _{F\times E} \} \cup
(P + \Pic ^0 (E)).$$
In particular $V^0(X,\Pic ^0(E),\omega _X  )$ does not generate
$\Pic ^0(E)$.

\proclaim{Corollary 2.5} Let $X$ be any variety of maximal Albanese dimension.
$\kappa (X)$ is invariant under smooth deformations. 
\endproclaim

\demo{Proof}
Let $\Delta $ be an open neighborhood of 
a point $0$ of a smooth projective curve, 
$\delta :\XX \lra \Delta$ be a smooth projective morphism
with connected smooth
fibers, $\XX _0:=\delta ^{-1}(0)$ be a closed fiber of $\delta$
such that $\XX _0=X$. Since $q(X)$ is deformation invariant,
for all $t\in \Delta$, $A_t:=\Alb (\XX _t)$ is an abelian variety of 
dimension $q(X)$. 

Let $P\in \Pic^0(X)$ such that $P^{\ot m}=\OO _X$,
$\PP$ be a section of $\Pic ^0(\XX /\Delta )$ such that
$P=P_0:=\PP |_{\XX _0}$,
$P_t:= \PP |_{\XX _t} \in \Pic ^0(A_t)$ satisfies $P_t^{\ot m}=\OO _{\XX _t}$.
Let 
$$\tilde {\XX _t} :=Spec  \left( \oplus _{i=0}^{m-1}P_t^{\ot i}\right)$$ 
be the corresponding \'etale cyclic cover of degree
$m$.
We have that the quantity
$$h^0(\tilde {\XX _t}, \omega _{\tilde {\XX _t}})=\sum_{i=0}^{m-1}
h^0(\omega _{{\XX_t}}\ot P_t^{\ot i})$$
is constant. However the functions $h^0( {\XX _t}, \omega _{{\XX_t}}\ot P_t^
{\ot i})$ are upper
semicontinuous in $t\in \Delta$, and hence also constant.

Let $\overline {T} ^i$ denote the translates
of the components $T^i$ of $V^0(X,\Pic ^0(A),\omega _X)$
through the origin. The subvarieties $T^i$ are determined by their
torsion points (cf. Theorem 1.3.1.d). 
In particular, recall that
$\Pic ^0(\XX _t)\cong H^1(\XX _t, \OO _{\XX _t}^*)/H_1(\XX _t,\ZZ)$
and the subvarieties $\overline {T}^i_t$ are determined by the corresponding
vector subspaces of $H_1(\XX _t,\QQ )$.
We remark that given two $\QQ$ vector subspaces $W_t^i\subset H_1(\XX _t,\QQ )$
continuous in the parameter $t$, then $\dim (W_t ^1\cap W_t ^2)$ and
$\dim (W_t ^1+ W_t ^2)$ are constant.

{From} the above discussion, it follows
that there exist subvarieties $T^i_t$ of
$\Pic ^0(A_t)$, which are smooth deformations of $T^i_0$ such that
$V^0(\XX _t,\Pic ^0(A_t),\omega _{\XX _t})=\cup T^i_t$.
Moreover, for any set of
indices $I$, the quantities
$$\dim \left( \bigcap _{i\in I} \overline {T} ^i_t\right) \ \ \ and\ \ \
\dim \left( \sum _{i\in I} \overline {T} ^i_t \right)  $$
are constant.
In particular the quantity
$$\dim \left( V^0(\XX _t,\Pic ^0(A_t), \omega_{ \XX_t}) \right) =\kappa (\XX _t)+q(\XX _t)
-\dim (\XX _t)$$
is constant and hence $\kappa (\XX _t)$ is also constant.
\hfill $\square$ \enddemo

\heading 3. Pluricanonical maps of
varieties of maximal Albanese dimension \endheading

We will keep the notation of the preceding sections.
In particular $X$ will be a smooth projective variety with maximal
Albanese dimension and 
$a:X\lra A$ will denote a generically finite morphism to an abelian variety.

\proclaim{Lemma 3.1} Let $E$ be an $a$-exceptional effective
divisor on $X$.
If $\OO _X (E)\ot P$ is effective, then $P=\OO _X$.
\endproclaim
\demo{Proof}
If $\dim (X)=1$, then $E=0$ and the assertion is clear.
If $\dim (X)>1$,
pick $H_A$ a sufficiently ample divisor on $A$ and let $H=a^*H_A$ be
the corresponding nef and big divisor on $X$. Choosing $H_A$ appropriately,
we may assume that $H$ is a smooth divisor in $X$. We may also assume that
$h^0(X,\OO _X(E-H) \ot P)=0$ for all $P\in \Pic ^0(X)$.
If $\dim (X)\geq 2$, from the exact sequence of sheaves

$$0\lra \OO _X(E-H) \ot P\lra \OO _X(E) \ot P\lra \OO _H(E) \ot P\lra 0,$$
we have an injection
$H^0(X,\OO _X(E) \ot P)\hookrightarrow  H^0(H,\OO _H(E) \ot P)$ for all
$P\in \Pic ^0(X)$. Similarly, there is an exact sequence of sheaves

$$0\lra \OO _X(-H) \ot P\lra \OO _X \ot P\lra \OO _H \ot P\lra 0.$$
Since $H$ is nef and big,
$$h^i(X,\OO _X(-H) \ot P)=h^{n-i}(X,\omega _X (H)\ot P^*)=0\ \ \ \ for \
all\ i<n.$$
Therefore, $h^0(X,\OO _X \ot P)=h^0(H,\OO _H
\ot P)=0$ for $P \ne \OO _X$ and $h^0(X,\OO _X )$ $=h^0(H,\OO _H)$ $=1$.
Clearly, for a general choice of $H$, we have that $H$ is of
maximal Albanese dimension, and $E|_H$ is an $a|_H$-exceptional divisor.
Repeating the above procedure, by successively intersecting 
appropriate divisors
pulled back from $A$, one obtains a curve $C\subset X$ such that

i) If $h^0(X, \OO _X (E)\ot P)>0$ then $h^0(C, \OO _C (E)\ot P)=
h^0(C, \OO _C \ot P)>0$,

ii) If $h^0(C, \OO _C \ot P)>0$ then $P=\OO _X$.

\noindent It follows that if $h^0(X, \OO _X (E)\ot P)>0$, then $P=\OO _X$.
\hfill $\square$ \enddemo

\proclaim{Lemma 3.2} Let $D$ be an irreducible reduced divisor on $X$
which is not $a:X\lra A$ exceptional, $H$ a Cartier divisor on $X$ which is
numerically
trivial on the general fiber of $X\lra S$. Then, $D$ is not contained
in the base locus of $|mK_X+H|$ for infinitely many values of $m$.
\endproclaim
\demo{Proof} \cite{Ko3, 17.6.1}.
\hfill $\square$ \enddemo
Let $|mK_X|=F_m+|M_m|$, where $Bs|M_m|$ contains no divisors.
We remark that, since $X$ is of maximal Albanese dimension,
we may assume that $K_X$ is effective.

\proclaim{Theorem 3.3} If $X$ is of general type, i.e. $\kappa (X)=\dim
(X)$.
Then for any integer $s\geq 3$, $F_s$ is
$a$-exceptional. $|5K_X|$ defines
a generically finite rational map and $|(5\dim (X)+1)K_X|$ defines
a birational map.
\endproclaim
\demo{Proof}
For any fixed divisor $D$ on $X$,
then there
exists an $m_0$ such that $|mK_X-D| \ne \emptyset$ for all $m \ge
m_0$. Let $R$ be an irreducible divisor of $X$ which is not
$a$-exceptional. By Lemma 3.2,
$R \not \subset Bs|mK_X-D|$ for infinitely many $m
\gg 0$.
\proclaim{\bf Step 1} $F_s$ is exceptional for all $s \ge 3$.
\endproclaim
It is easy to see that $F_s$ being $a$-exceptional is a condition
which is independent of the
particular birational model of $X$ under consideration.
$F_1\subset R_a $ , the ramification divisor of $a: X \lra A$, contains at
most finitely many non
$a$-exceptional components, which we denote by $R_i$.
Fix positive integers $m_i$ such that $R_i\not \subset Bs |m_iK_X|$,
then $R_i\not \subset Bs |\lambda m_iK_X|$ for any integer
$\lambda >0$. It follows that for $m_0=\prod m_i$, $F_{m_0}$ is
$a$-exceptional.

Next, fix an ample divisor $H$ on $A$, then $a^*H$ is nef and big on $X$.
Let  $K_X$ be a canonical divisor.
Let $r$ be the multiplicity of $K_X$ at $R$. We distinguish two cases.

\noindent {\bf Case 1.} $r= 1$.

By Lemma 3.2, there exists a positive integer $t$ such that the base locus
of
$|tK_X -a^*H|$ doesn't contain $R$. Let $B$ be an general element
of $|(s-2)(tK_X -a^*H)|$. By replacing $X$ with an appropriate birational
model, we may assume that $R$ is smooth and $B+R$ has normal crossing support
and $B$ does not contain $R$. Define $$M:= \OO_X\left( (s-2)K_X- \lfloor
\frac{B}{t} \rfloor \right)
\equiv \left( \frac{s-2}{t} \right) a^*H+ \{ \frac{B}{t}
\}.$$
Consider the exact sequence
$$ 0 \to \omega _X\ot M \to \omega _X \ot M(R) \to \omega _R \ot M \to 0. $$
Since $a^*H$ is nef and big and $\{ \frac{B}{t} \}$ is klt,
one sees that $H^1(X,\ox \ot M)=0$.
The divisor $R$
is not $a$-exceptional, so $K_R$ is effective and $(a^*H)|_R$ is nef and
big. $R$ is not contained in the support of $B$, so the $\QQ$ divisor
$\{ \frac{1}{t} (B|_R )\}$ is also klt. By Theorem
1.1.1.e, $H^0(R,\omega_R \ot M) \ne 0$. Therefore, there is a divisor in
$|K_X
+M+R|$ not containing $R$. This gives  a divisor in $|sK_X|=
|K_X+M+R+ (K_X -R)+\lfloor \frac{B}{t} \rfloor|$ not containing
$R$.

\noindent {\bf Case 2.} $ r\ge 2$.

There exists positive integers $t,\ m_0$ such that
$R \not \subset Bs|tK_X -a^*H|$, and $R \not \subset
Bs|m_0K_X|$. For any integer $s\geq 2$ let $K':= (s-1)K_X -R$.
Consider the following linear series
$$ |mK'-m_0a^*H|= |\big((s-1)m-\frac{m}{r}-m_0t \big)
K_X+\frac{m}{r}(K_X-rR)+m_0(tK_X-a^*H)|.$$
It follows that for $m$ divisible by $m_0r$,
$R$ is not in the base locus of $|mK'-m_0a^*H|$.
Choose a general $B \in |mK'-m_0a^*H|$ and define

$$M:= \OO _X \left( K'- \lfloor \frac{B}{m} \rfloor \right)
\equiv \frac{m_0}{m} a^*H+
\{ \frac{B}{m} \}.$$
An argument  similar to the one in the previous case again
shows that $R$ is not in the base
locus of $|sK_X|$.
\proclaim {\bf Step 2} $|5K_X|$ defines a generically finite map.
\endproclaim
Replacing $X$ by an appropriate birational model, we may assume that
$|M_3|$ is base point free. Let $D$ be a general member of $M_3$ and
$\Delta$ be the image of $D$ in $A$. If $\Delta$ is not of general type,
then there exists an ample line bundle $H$,
a semipositive line bundle $L$ and a positive integer $s$
such that
$\Delta \cdot L^{s}\cdot H^{n-s-1}=0$ and $L^{s}\cdot H^{n-s}>0$.
It follows that also $D\cdot (a^*H)^{n-s-1}\cdot (a^*{L})^s=0$.

Since $X$ is of general type, there exists a rational
number $\epsilon >0$ such
that $K_X-\epsilon (a^*H)$ is an effective $\QQ$-divisor.
Therefore,
$$K_X\cdot (a^*H)^{n-s-1}\cdot (a^*{L})^s\geq \epsilon (a^*H)^{n-s}\cdot
(a^*{L})^s>0.$$
Since $K_X-D\equiv F_3$ is
$a$-exceptional, it follows that
$$D\cdot (a^*H)^{n-s-1}\cdot (a^*{L})^s=K_X\cdot (a^*H)^{n-s-1}
\cdot (a^*{L})^s>0.$$
This is a contradiction. We may therefore assume that $\Delta$ and hence
$D$
are of general type. It is well known that for an appropriate
desingularization
$\Delta '$ of $\Delta$, the linear series $|K_{\Delta '}|$ defines a
generically finite rational map.
Therefore, we may assume that $|K_D|$ defines a
generically finite rational map. Let $B\in |mK_X-a^*H|$,
and $M:=\OO _X(K_X-\lfloor \frac{B}{m}\rfloor)\equiv 
\frac{a^*H}{m}+\{\frac{B}{m}\}$.
We may assume that $B$ has normal crossings support and does not contain
$D$.
Consider the exact sequence of sheaves
$$0\lra \ox \ot M\lra \ox \ot M(D)\lra \omega _D\ot M\lra 0,$$
which is also exact on global sections. 

\proclaim{Claim} There exist $m,B$ such that $M$ is effective.
\endproclaim
Fix $m_1$ such that $|m_1K_X-a^*p^* H_S|$ is non empty. Let $B_1\in 
|m_1K_X-a^*p^* H_S|$. Let $m:=m_1+m_2$ and $B:=B_1+m_2K_X\in|m 
K_X-a^*p^* H_S|$. Let $\Gamma $ be any prime divisor, and $b_1,k$ be the 
multiplicities of $B_1,K_X$ along $\Gamma$. For $m_2\gg 0$, we have that
$$\lfloor \frac{b_1+km_2}{m}\rfloor =\lfloor \frac{b_1}{m_1+m_2}+k
\frac {m_2}{m_1+m_2}\rfloor\leq k.$$
Since there are only finitely many components of $B+K_X$, for $m_2\gg 0$, 
we have $\lfloor \frac{B}{m} \rfloor \prec K_X$, and the claim follows.

It is easy to see that $M|_D$ is also effective.
Therefore,
$|K_X+D+M|$ restricted to $D$ defines a generically finite rational map.
If $|K_X+D+M|$ does not define a generically finite rational map, then
the closure of the image of $X$ 
which we denote by $Y\subset \PPP =\PPP (H^0(X,
\ox\ot M(D)
))$, must be dominated by $D$. This is however impossible
since by Theorem 1.1.1.e $h^0(X,\ox  \ot M)>0$, and hence $|K_X+D+M|$
contains non-trivial sections vanishing on $D$, i.e. there is a hyperplane
section of $\PPP$ containing the image of $D$ but not containing $Y$.
Finally, the assertion follows 
from the inclusion of sheaves
$$\ox \ot M(D) \hookrightarrow  \ox^{\ot 5}.$$

By \cite{Ko1, Theorem 8.1},
$|(5\dim (X)+1)K_X|$
defines a birational map.
\hfill
$\square$
\enddemo

\proclaim{Theorem 3.4} $|6K_X|$ defines
rational map with image of dimension $\kappa (X)$,
and $|(6 \kappa(X)+2)K_X|$ defines
the stable canonical map.
\endproclaim
\demo{Proof}
Let $H_A$, $H_S$ be the pull backs of sufficiently ample divisors on
$A$, $S$ respectively. For any $P\in \Pic ^0(X)$, let
$$|mK_X+P|=|M_{m,P}|+F_{m,P},$$
where $Bs|M_{m,P}|$ contains no divisor.

\proclaim {Step 1} Let $S\lra S'$ be any surjective map of abelian varieties
such that $y:=\dim (S)>y':=\dim (S')$ and $s':X \to S'$ be the induced map.
Then, for any divisor $D\in |M_{4}|$,
$D$ is not $s'$-vertical.
\endproclaim
By Lemma 1.3.4, the linear series $|2K_X+P|$ is nonempty for all
$P\in \Pic ^0(S)$.
Let  $Z$ be an irreducible component of $H^{y'}_{S'}$, ie a
general fiber of $s':X \to S'$.
For all $P\in \Pic ^0(S)$ one has a map of linear series
$$|2K_X+P|\times |2K_X-P|\lra |4K_X|.$$
It follows that
$$F_4\subset \bigcap _{P\in \Pic ^0(S)} (F_{2,P}+F_{2,-P}).$$

\proclaim {Claim 1} $\lceil \frac{F_4}{2}\rceil \prec F_{2,P}$
for general $P\in \Pic ^0(S)$.
\endproclaim
Let $V\subset \Pic ^0(S)$ be an open set such that
$h^0(X,\omega _X ^{\ot 2}\ot P)$ 
is constant. Let $R$ be any component of $F_4$ of multiplicity
$r$, then $\lceil \frac{r}{2}\rceil R \prec F_{2,P}$ for all
$P\in U\subset \Pic ^0(S)$
where $U$ is a Zariski dense subset such that
$U\cup -U=\Pic ^0(S)$. 
Since the condition $\lceil \frac{r}{2}\rceil R \prec Bs
|2K_X+P|$
is Zariski closed in $V$, it follows that $\OO _X(2K_X-\lceil 
\frac{r}{2}\rceil R )
\ot P$
is effective for all $P\in V\subset \Pic^0(S)$. Since this holds for all
components of $F_4$, the claim follows.

\proclaim {Claim 2} $ker ( \Pic^0(S) \to \Pic^0(Z) )$ is a proper closed
subvariety
of $\Pic^0(S)$.
\endproclaim
Let $V$ be the image of $Z$ in $S$.
$V$ is a translate of $ker(S \to S')$. There exists infinitely many 
$P \in \Pic^0(S)$
such that $P|_V \ne \OO_V$. Hence pulling back to $Z$, we have
$P|_Z \ne \OO_Z$.

\proclaim {Claim 3} $(2K_X-\lceil \frac{F_4}{2}\rceil )|_Z$ is not $a|_Z$-exceptional.
\endproclaim

By Claim 1, $(2K_X-\lceil \frac{F_4}{2} \rceil +P)|_Z$ is
effective and for general $P \in \Pic^0(S)$.
By Claim 2, we may assume that $P|_Z \ne \OO_Z$ for general $P \in \Pic^0(S)$.
Since $a|_Z:Z\lra a(Z)$ is generically finite, the claim follows by Lemma 3.1.

\bigskip
Claim 3 implies that
$2K_X-\lceil \frac{F_4}{2}\rceil$ is not $s'$-vertical
(otherwise $(2K_X-\lceil \frac{F_4}{2}\rceil )|_Z$ $=\OO _Z$).
{From} the inclusion of linear series
$$|2K_X-\lceil \frac{F_4}{2}\rceil |\times |2K_X-\lceil \frac{F_4}{2}\rceil |
@>{+2\lceil \frac{F_4}{2}\rceil -F_4}>>|4K_X-F_4|$$
it follows that $4K_X-F_4=M_{4}$ is also not $s'$-vertical.

\proclaim {Step 2} Let $D$ be a general member of $M_{4}$, and $\Delta$ be
its image in $S$. Then $\Delta$ is of general type.
\endproclaim
If $\Delta$ is not of general type, then
$\Delta$ is vertical for an appropriate projection $S\lra S'$.
By Step 1 this is impossible.

\proclaim {Step 3} $|6K_X|$ defines a
rational map whose image is of dimension $\kappa(X)$.
\endproclaim
Let $\tilde{\Delta}$ be an appropriate desingularization of $\Delta$.
The linear series $|K_{\tilde{\Delta}}|$ defines a 
generically finite rational map.
Replacing $X$ by an appropriate birational model, we may assume that $|M_4|$
is free and hence $D$ is a smooth subvariety that maps onto $\tilde{\Delta}$.

Fix an ample divisor $H_S$ on $S$.
Let $B\in |mK_X-a^*p^*H_S|$ be a general member.
Replacing $X$ by an appropriate birational model, we may assume that
$B$ has normal crossings support.
Define $$M:=\OO _X\left( K_X-\lfloor \frac{B}{m} \rfloor \right)
\equiv \frac{a^*p^*H_S}{m}+\{ \frac{B}{m}  \}.$$ 
As in the proof of the previous theorem, we may assume
that $M$, $M|_D$ are effective,
and $B$ does not contain $D$. Consider the exact sequence:

$$0\lra \ox \ot M\lra \ox \ot M(D)\lra \omega _{D}\ot M\lra 0.$$
By \cite{Ko3,Theorem 10.19},
this is exact on global sections.
Sections of $\omega_{D} \ot M$ lift to sections of $\omega _X \ot M(D)$.
Since $\OO _{ D}(K_{  D/\tilde{\Delta}})\ot M$ is effective.
It follows that $\omega _{D}\ot M$ also defines a rational map
with image of dimension at least $\kappa(X)-1=\dim (W) -1$.
By Theorem 1.1.1.e, $|K_X+M|$ is non empty. An argument similar 
to the one in the proof of Theorem 3.3
shows that  $|K_X+M+D|$ defines a
rational map
with image of dimension at least $\kappa(X)$.
The assertion follows from the inclusion of sheaves
$$\ox \ot M(D) \hookrightarrow  \ox^{\ot 6}.$$
\proclaim{Step 4} $|(6\kappa (X)+2)K_X|$ is the stable
canonical map.
\endproclaim
We may assume that $|6K_X|=|M_6|+F_6$ and $M_6$ is free.
Let $M$ be defined as in step 3. $\varphi _6$ factors as $X\lra Y\lra Y'
:=\varphi _6(X)$.
Pick $D_1,...,D_{\kappa}$
general sections of $M_6$, with $\kappa =\kappa (X)$.
Let $X_i:=D_1\cap...\cap D_i$. Then $X_{\kappa }=D_1\cap ...\cap
D_{\kappa}$ is the union of $\deg (Y\lra Y')$ fibers of $X\lra Y$
which we denote by $F_t$.
We must show that sections of $|(6\kappa (X)+2)K_X|$ separate these fibers.
Let $g_i:X_i@>>> \bar {X}_i$ be the maps induced by $X@>>> Y'$.
By \cite{Ko3, Theorem 10.19}, we have that
$$H^1( {X}_i,\omega _{X_i}\ot M\ot M_6^{\ot \kappa -i-1}) \lra
H^1( {X}_i,\omega _{X_i}\ot M\ot M_6^{\ot \kappa -i})$$
 is injective for all  $i\leq \kappa -1$.
Therefore, the exact sequences
$$0 \to \omega _{X_i}\ot M\ot M_6^{\ot \kappa -i-1}\to \omega _{X_i}\ot
M\ot M_6^{\ot \kappa -i}\to 
\omega _{X_{i+1}}\ot M\ot M_6^{\ot \kappa -i-1}  \to 0,$$
are exact on global sections.
It follows that $$H^0(X, \omega _X \ot M\ot M_6^{\ot \kappa(X)})\lra 
H^0(X_{\kappa }, \omega _{X_{\kappa }}\ot M)$$ is surjective.
The assertion now follows 
since $\omega _X \ot M\ot M_6^{\ot \kappa(X)}$ 
is a subsheaf of $\ox^{ \ot ( 6 \kappa(X)+2)}$,
$H^0(X_{\kappa }, \omega _{X_{\kappa }}\ot M)$ 
$=\oplus H^0(F_t,\omega _{F_t}\ot
M )$, and by Theorem 1.1.1.e, $H^0(F_t,\omega _{F_t}\ot
M ) \ne 0$.
\hfill $\square$
\enddemo

\enddocument
\end